\documentclass[10 pt, reqno]{amsart}
\usepackage {amsfonts}
\usepackage{amsthm}
\usepackage{amssymb}
\usepackage{latexsym}
\usepackage{amsmath}
\usepackage{mathrsfs}
\usepackage{comment}

\theoremstyle{plain}
\newtheorem{tw}{Theorem}[section]

\newtheorem {lem} [tw]{Lemma}
\newtheorem {prop}[tw] {Proposition}

\theoremstyle{definition}
\newtheorem*{deft}{Definition}
\theoremstyle{remark}
\newtheorem*{rem}{Remark}
\newtheorem*{rems}{Remarks}
\newcommand{\bc} {\mathbb C}
\newcommand{\bn}{\mathbb N}
\newcommand{\br}{\mathbb R}

\newcommand{\bialg} {\mathsf{B}}
\newcommand{\Al} {\mathsf{A}}

\newcommand{\clg} {\mathsf{C}}

\newcommand{\semiGp} {H}

\newcommand{\mlg} {\mathsf{M}}
\newcommand {\id} {{\textup{id}}}

\newcommand{\Ind}{\mathcal{I}}

\newcommand{\multib}{M(\bialg)}

\newcommand{\tu}{\textup}

\newcommand{\hil}{\mathsf{h}}

\newcommand{\strfunb}{M(\bialg)_{\beta}^*}
\newcommand{\gen}{Z}

\newcommand{\Com}{\Delta}
\newcommand{\Cou}{\epsilon}

\newcommand {\Ker} {{\textup{Ker}}}
\newcommand {\Dom} {{\textup{Dom}}}

\newcommand{\Ran}{\textup{Ran}}

\newenvironment{rlist}
{

\begin{enumerate}}
{\end{enumerate}}

\newenvironment{alist}
{

\begin{enumerate}}
{\end{enumerate}}

\newcommand{\la}{\lambda}

\newcommand{\La}{\Lambda}

\newcommand{\ot}{\otimes}

\newcommand{\otol}{\overline{\ot}}

\newcommand{\idAl}{\textup{id}_{\Al}}
\newcommand{\idb}{\textup{id}_{\bialg}}
\newcommand{\idm}{\textup{id}_{\mlg}}
\newcommand{\wt}{\widetilde}

\newcommand{\cb}{{\text{\tu{cb}}}}
\newcommand{\bd}{{\text{\tu{b}}}}

\newcommand{\tilded}[1]{{#1}_{\kern-1pt\raisebox{-1.25ex}{\mbox{\,$\tilde{\ }$}}}}
\newcommand{\otM}{\otimes_{\mathrm{M}}}
\newcommand{\Mot}{{\ }_{\mathrm{M}}\!\!\otimes}

\numberwithin{equation}{section}


\begin{document}
\author{J.\ Martin Lindsay and Adam G.\ Skalski}
\address{Department of Mathematics and Statistics,
Lancaster University, Lancaster LA1 4YF, U.K.}
\email{j.m.lindsay@lancaster.ac.uk}
\email{a.skalski@lancaster.ac.uk}

\title{\bf Convolution semigroups of states}

\begin{abstract}
Convolution semigroups of states on a quantum group
form the natural noncommutative analogue of
convolution semigroups of probability measures on a locally compact group.
Here we initiate a theory of
weakly continuous convolution semigroups of functionals on a $C^*$-bialgebra,
the noncommutative counterpart of a locally compact semigroup.
On locally compact quantum groups we obtain a
bijective correspondence between such convolution semigroups and
a class of $C_0$-semigroups of maps which we characterise.
On $C^*$-bialgebras of discrete type we show that all weakly continuous
convolution semigroups of states are automatically norm-continuous.
As an application we deduce a known characterisation of
continuous conditionally positive-definite Hermitian functions
on a compact group.
\end{abstract}

\maketitle
\section*{Introduction}

Convolution semigroups of probability measures on locally compact groups,
and semigroups,
play an important role in both probability theory and functional analysis.
They underpin all approaches to stochastic processes with
independent identically distributed increments, that is L\'evy processes
(\cite{Bertoin}).
The analysis of convolution semigroups is also an independent area of research,
well documented in the monograph~\cite{Heyer}.
Given the recent development of a satisfactory theory of locally compact
quantum groups (\cite{KV})
it is natural to investigate convolution semigroups of
states on such objects---states being
the noncommutative counterpart of probability measures.

We recently introduced the notion of quantum L\'evy process on a
locally compact quantum semigroup, in other words a $C^*$-bialgebra
(\cite{qscc3}).
As in the classical theory of L\'evy processes,
and the earlier noncommutative theory of L\'evy processes on
purely algebraic
quantum semigroups (\cite{Schurmann}),
the concept of a convolution semigroup of states plays a fundamental role
here too.
In our papers the strongest results are obtained under the assumption that
the Markov semigroup of the process is norm-continuous.
This class of processes may be viewed as
a quantum analogue of the class of compound Poisson processes.

The present paper forms the first step
in analysing general quantum L\'evy processes.
We initiate a theory of weakly continuous convolution semigroups of functionals
on locally compact quantum semigroups.
In particular we prove that on $C^*$-bialgebras enjoying
a `residual vanishing at infinity' property (which
is satisfied by all locally compact quantum groups)
every weakly continuous convolution semigroup of functionals
induces a $C_0$-semigroup on the algebra.
From a probabilistic viewpoint this amounts to a Feller property (\cite{Jacob})).
We give several characterisations of the resulting class of $C_0$-semigroups.
We also prove that on $C^*$-bialgebras of discrete type
weak continuity for a convolution semigroup of states implies norm-continuity.
Here slice maps for non-normal functionals,
in the sense of Tomiyama (\cite{Tomiyama}), play a crucial role.
Our results show that all quantum L\'evy processes on a discrete quantum semigroup
are automatically of the type thoroughly studied in~\cite{qscc3}.
As an application we prove a result of Guichardet, characterising the
continuous conditionally positive-definite Hermitian functions on a compact group
(\cite{Guichardet}).

\section{Preliminaries} \label{prelim}

Given $C^*$-algebras $\Al$ and $\clg$,
their spatial/minimal (respectively algebraic) tensor product
will be denoted $\Al \ot \clg$ (resp. $\Al\odot\clg$) and
the multiplier algebra of $\Al$ will be denoted by $M(\Al)$.
The multiplier algebra is equipped with the strict topology,
with respect to which the unit ball of $\Al$ is dense in that of $M(\Al)$.
A linear map $T$ from $\Al$ (or $M(\Al)$) to $M(\clg)$ is \emph{strict} if
it is bounded and strictly continuous on bounded subsets of $\Al$
(respectively $M(\Al)$).
Standard examples of strict maps include
$^*$-homomorphisms which are nondegenerate (that is $T(\Al)\clg$ is total in $\clg$),
continuous linear functionals on $\Al$ and, more generally,
slice maps of the type $\mu\ot \id_{\clg}: \Al \ot \clg \to \clg \subset M(\clg)$,
where $\mu \in \Al^*$.
The collection of strict linear functionals on $M(\Al)$
forms a closed subspace of $M(\Al)^*$ which we denote by $M(\Al)^*_\beta$.
An important property of strict maps $T: \Al \to M(\clg)$
is that they possess a unique strict extension $\wt{T}: M(\Al) \to M(\clg)$;
the extension procedure does not change the norm
(or cb-norm if $T$ is completely bounded) and preserves (complete) positivity.
In particular, strict extension defines an isometric isomorphism $\Al^* \to M(\Al)^*_\beta$
mapping states to states.
A strict map $T$ is called \emph{preunital} if its strict extension is unital: $\wt{T}(1)=1$.
Thus all states are preunital, as are nondegenerate *-homomorphisms.
The existence and uniqueness of extensions allows \emph{composition}
of strict maps $T:\Al_1 \to M(\Al_2)$ and $S:\Al_2 \to M(\Al_3)$,
with the resulting map $\Al_1 \to M(\Al_3)$ being strict.
Finally, for strict completely positive maps $T_i:\Al_i\to\clg_i$ ($i=1,2$),
their tensor product map
$T_1\ot T_2: \Al_1\ot\Al_2 \to M(\clg_1)\ot M(\clg_2)\subset M(\clg_1\ot\clg_2)$, is strict.
For more information on multiplier algebras, strict maps, their tensor products and
extensions, we refer to~\cite{Lance},~\cite{Johan} and~\cite{qscc3}.

By a \emph{semigroup on} $\Al$, or $M(\Al)$,
we mean a one-parameter family of bounded operators $(P_t)_{t\geq 0}$
satisfying $P_{s+t}=P_s P_t$
and $P_0 = \idAl$
(no continuity in $t$ is assumed a priori).
We call it a \emph{strict/completely bounded/completely positive} semigroup
if each map $P_t$ has that property.
By elementary semigroup theory (\cite{Davies}, Proposition 1.23),
$(P_t)_{t\geq 0}$ is a $C_0$-semigroup, in other words it is strongly continuous in $t$,
if it is weakly continuous, that is
$t\mapsto \mu(P_t(a))$ is continuous for each $a\in\Al$ and $\mu\in\Al^*$,
at $t=0$.

We shall need two general results,
in Sections~\ref{on Cstar Bialgebras} and~\ref{disc type}
respectively.
The first concerns
the preservation of pointwise convergence under strict extension.

\begin{prop} \label{contstat}
Let $(\omega_\lambda)_{\lambda \in \Lambda}$ be
a net of states on a $C^*$-algebra $\Al$
converging pointwise to a state $\omega$ on $\Al$.
Then $(\wt{\omega}_\lambda)_{\lambda \in \Lambda}$ converges pointwise to $\wt{\omega}$.
\end{prop}
\begin{proof}
Let $m \in M(\Al)$ and let $(e_i)_{i \in \Ind}$ be an approximate unit for $\Al$.
Then, for each $i\in\Ind$,
\[
\wt{\omega}_\lambda ((1- e_i)^2) =
1 + \omega_\lambda (e_i^2-2 e_i)
\stackrel{\lambda }{\longrightarrow} 1 + \omega (e_i^2-2 e_i) =
\wt{\omega}((1-e_i)^2)
\]
and so, by the Cauchy-Schwarz inequality for states,
\begin{align*}
|(\wt{\omega}_\lambda - \wt{\omega})(m)|
&\leq
|\wt{\omega}_\lambda (m(1-e_i))| + |\wt{\omega} (m(1-e_i))| +
|(\omega_\lambda - \omega)(me_i)|
\\&\leq
\|m\|(\wt{\omega}_\lambda ((1-e_i)^2))^{1/2} +
\|m\| (\wt{\omega} ((1-e_i)^2))^{1/2}+
|(\omega_\lambda - \omega)(me_i)|
\\&
\stackrel{\lambda }{\longrightarrow}
2 \|m\| \wt{\omega}\big((1-e_i)^2\big)^{1/2}.
\end{align*}
Since $\big((1-e_i)^2\big)_{i\in\Ind}$ converges strictly to $0$, the result follows.
\end{proof}

\begin{rem}
Using Kadison-Schwarz, instead of Cauchy-Schwarz, the above proof remains valid for
strict completely positive preunital maps into the multiplier algebra of another
$C^*$-algebra $(T_\lambda: \Al\to M(\clg))_{\lambda\in\Lambda}$
converging pointwise to another such map $T$.
\end{rem}

The following observation is due to Tomiyama (\cite{Tomiyama});
we include its elementary proof for the convenience of the reader.
The notation $\otol$ is used for ultraweak tensor products of
von Neumann algebras
and ultraweakly continuous slice maps (\cite{KaR}).

\begin{prop}\label{new}
Let $\mlg = \mlg_1\otol\mlg_2$,
for von Neumann algebras $\mlg_1$ and $\mlg_2$,
and let $\nu\in (\mlg_1)^*$.
Then there is a unique map $\nu\otM\id_2: \mlg\to\mlg_2$
satisfying
\begin{equation}\label{U}
\varphi\circ (\nu\otM\id_2) =
\nu\circ (\id_1\, \otol\, \varphi),
\quad \varphi\in (\mlg_2)_*,
\end{equation}
where $\id_i$ denotes the identity map on $\mlg_i$ \tu{(}$i=1,2$\tu{)}.
Moreover,
$\nu\otM \id_2$ is a bounded operator of norm $\|\nu\|$.
\end{prop}


\begin{proof}
Let $x\in\mlg$. Then, for $\varphi\in (\mlg_2)_*$,
\[
\big| \nu \big((\id_1\, \otol\, \varphi )(x)\big)\big| \leq
\|\nu\| \big\|\id_1\, \otol\, \varphi \big\| \| x \| =
\|\nu\| \|\varphi \| \| x \|
\]
so the prescription $\varphi \mapsto \nu ((\id_1\, \otol\, \varphi )(x))$
defines a bounded linear functional on $(\mlg_2)_*$
of norm at most $\| \nu\| \| x \|$.
Therefore, invoking the canonical identification of
$((\mlg_2)_*)^*$ with $\mlg_2$,
there is a bounded operator $\nu\, \otM\, \id_2 : \mlg \to \mlg_2$,
of norm at most $\|\nu\|$, satisfying~\eqref{U}.
Since $\nu\, \otM\, \id_2$ clearly extends the map
$\nu\, \ot\, \id_2 : \mlg_1 \ot \mlg_2 \to \mlg_2$, which has norm $\|\nu \|$,
it follows that $\|\nu\, \otM\, \id_2\| = \|\nu\|$.
\end{proof}
\begin{rems}
By uniqueness, and the fact that
$(\mlg_2)_* = \{ \phi|_{\mlg_2}: \phi \in B(\hil)_* \}$
where $\hil$ is the Hilbert space on which $\mlg_2$ acts, it follows that
$\nu \otM \id_2 = \nu \otM \id_{B(\hil_2)}|_{\mlg}$.
In particular, by Remark 1.4 of~\cite{SteveMartin},
$\nu \otM \id_2$ is completely bounded with
cb-norm equal to $\|\nu\|_{\cb} = \|\nu\|$;
moreover the map $\nu \otM \id_2$ is completely positive if
the functional $\nu$ is positive.
In fact, the above result still holds if $\nu$ is replaced by a
completely bounded map into another von Neumann algebra $\mlg_3$,
with~\eqref{U} now reading
$(\id_3\, \otol\, \varphi) (\nu\otM\id_2) =
\nu\circ (\id_1\, \otol\, \varphi)$ (\cite{SteveMartin}).

If $\nu$ is ultraweakly continuous then
$\nu \otM \id_2$ equals $\nu\, \otol\, \id_2$.
Naturally there are also bounded operators
$\id_1 \Mot \mu : \mlg\to\mlg_1$, for each $\mu\in (\mlg_2)^*$,
and the validity of the commutation relation
\begin{equation}\label{comm}
\nu_1 \circ (\id_1 \Mot \nu_2) = \nu_2 \circ (\nu_1 \otM \id_2)
\end{equation}
for functionals $\nu_i \in (\mlg_i)^*$, $i=1,2$, naturally arises.
Under the assumptions that
$\mlg_2$ is infinite dimensional and
the Hilbert space on which $\mlg_1$ acts is separable,
Tomiyama showed that~\eqref{comm} holds for all
$\nu_2 \in (\mlg_2)^*$ if and only if
$\nu_1$ is ultraweakly continuous, in other words $\nu_1\in(\mlg_1)_*$
(\cite{Tomiyama}, Theorem 5.1).
There is a corresponding result for completely bounded maps
(\cite{Neufang}, Theorem 5.4).
\end{rems}

\section{Convolution semigroups on $C^*$-bialgebras} \label{on Cstar Bialgebras}

\noindent
We first recall the definition.

\begin{deft}
A $C^*$-\emph{bialgebra}
is a $C^*$-algebra $\bialg$ equipped with
a nondegenerate $*$-homomorphism $\Com:\bialg \to M(\bialg \ot \bialg)$, called the \emph{coproduct}, and
a character $\Cou: \bialg \to \bc$, called the \emph{counit}, satisfying
the coassociativity and counital properties:
\[
(\idb \ot \Com) \Com =  (\Com \ot \idb ) \Com
\ \text{ and } \
(\idb \ot \Cou) \Com = \idb = (\Cou \ot \idb ) \Com.
\]
A useful neutral expression for the first two maps is $\Com\!^{(2)}$. Thus $\Com\!^{(2)}: \bialg \to M(\bialg \ot \bialg \ot \bialg)$.
$C^*$-bialgebras are also called \emph{locally compact quantum semigroups};
unital $C^*$-bialgebras are called \emph{compact quantum semigroups}, or
$C^*$-\emph{bialgebras of compact type}.
\end{deft}

\emph{Fix now, and for the rest of the paper}, a (counital) $C^*$-bialgebra $\bialg$.
For $\la, \mu \in \bialg^*$, their \emph{convolution} is defined
as the following composition of strict maps:
\[
\la \star \mu = (\la \ot \mu) \Com.
\]
Convolution may also be viewed as a product on $\strfunb$.

\begin{prop} \label{banisom}
Both $(\bialg^*, \star)$ and $(\strfunb, \star)$ are unital Banach algebras, with
respective identities $\Cou$ and $\wt{\Cou}$, and
strict extension defines a unital isometric algebra isomorphism
from the former to the latter.
\end{prop}
\begin{proof}
Straightforward.
\end{proof}
\begin{deft}
A \emph{convolution semigroup of functionals} on $\bialg$ is
a family of functionals $(\lambda_t)_{t\geq 0}$  in $\bialg^*$ satisfying
\[
\la_0 =\Cou \ \text{ and } \ \la_{s+t}=\la_s\star\la_t,\quad s,t \geq 0;
\]
It is \emph{weakly continuous} if
\[
\lim_{t \to 0^+} \la_t(a) = \Cou(a), \quad a \in \bialg
\]
and \emph{norm-continuous} if
\[
\lim_{t\to 0^+} \|\la_t - \Cou\| = 0.
\]
\end{deft}

Restricting to semigroups of states,
Propositions~\ref{contstat} and~\ref{banisom} imply the following.

\begin{prop}
\label{2.2*}
Strict extension defines a one-to-one correspondence between
weakly continuous convolution semigroups of states on $\bialg$ and
weakly continuous convolution semigroups of strict states on $M(\bialg)$.
\end{prop}


The norm-continuous case is summarised in the next result. In brief,
norm-continuous convolution semigroups have bounded generators, from which the semigroup may be recovered by exponentiation. This also operates at the multiplier algebra level.

\begin{prop}
\label{char}
Let $(\la_t)_{t\geq 0}$ be a norm-continuous convolution semigroup of functionals on $\bialg$.
Then there is a unique functional $\gamma \in \bialg^*$ such that
\[
\lim_{t \to 0^+} \big\|\frac{1}{t}(\la_t - \Cou)-\gamma \big\|=0.
\]
Moreover,
\[
\la_t = \sum_{n=0}^{\infty}\frac{1}{n!}(t\gamma)^{\star n}, \;\;\;\;  \wt{\la}_t = \sum_{n=0}^{\infty}
\frac{1}{n!}(t\wt{\gamma})^{\star n}, \quad t\geq 0,
\]
and
\[
\lim_{t \to 0^+} \big\| \frac{1}{t}(\wt{\la}_t - \wt{\Cou})- \wt{\gamma} \big\|=0.
\]
\end{prop}
The functional $\gamma$ is called
the \emph{generating functional} of $(\la_t)_{t \geq 0}$.
\begin{proof}
In view of Proposition~\ref{banisom}, this follows from elementary properties of norm-continuous semigroups in a unital Banach
algebra.
\end{proof}

The class of norm-continuous convolution semigroups of states
is analogous to the class of classical convolution semigroups of
\emph{compound Poisson type}, with `initial measure' $\gamma + \Cou$.
More justification for this terminology is given in~\cite{Uwe} and~\cite{qscc2}.
The class of functionals $\gamma\in\bialg^*$ which generate a
convolution semigroup of states on $\bialg$ is characterised as follows
(\cite{qscc3}):
$\gamma$ is *-preserving,
\emph{conditionally positive}
in the sense that $\gamma (a)\geq 0$ for $a \in \bialg_+\cap \Ker\, \Cou$
and satisfies $\wt{\gamma}(1)=0$.
This is a form of Sch\"onberg correspondence;
see~\cite{Schurmann} for its algebraic (i.e.\ non-analytic) counterpart.
For further classical motivations we refer to~\cite{Heyer}.

To each functional $\mu \in \bialg^*$ there are associated two completely bounded strict maps $\bialg \to M(\bialg)$: 
\begin{equation}\label{L and R}
L_{\mu} = (\mu \ot \idb) \Com \ \text{ and } \
R_{\mu} = (\idb \ot \mu) \Com,
\end{equation}
and their strict extensions $\widetilde{L}_\mu, \widetilde{R}_\mu: M(\bialg)\to M(\bialg)$. Note that the original functional may
be recovered from either of these maps: $\Cou\circ L_\mu = \mu = \Cou\circ R_\mu$; from which it follows that $\|L_\mu\|_\cb =
\|L_\mu\| = \|\mu\| = \|R_\mu\| = \|R_\mu\|_\cb$. If the map $L_\mu$ is positive then the functional $\mu$ is positive, so
$L_\mu$ is actually completely positive;
also $\widetilde{L}_\mu$ is unital if and only if $\widetilde{\mu}$ is.
Thus $\widetilde{L}_\mu$ is completely positive and unital
if and only if $\mu$ is a state.
The same goes for the maps $R_\mu$ and $\widetilde{R}_\mu$.
The strict extensions of the $L$-map $\mu \mapsto L_\mu$ and
$R$-map $\mu \mapsto R_\mu$ define
completely isometric unital algebra morphisms
between the Banach algebras $(\bialg^*,\star)$
and $CB(M(\bialg))$, sharing the same left inverse: $\wt{T} \mapsto \Cou\circ T$, where $T := \wt{T}|_\bialg$.

For any convolution semigroup of functionals $(\la_t)_{t \geq 0}$ on $\bialg$, $(\wt{P}_t:=\wt{R}_{\la_t})_{t\geq 0}$ therefore
defines a strict semigroup of completely bounded maps on $M(\bialg)$, which we call the \emph{associated semigroup on the
multiplier algebra}. It determines the original convolution semigroup via the identity:
\begin{equation}\label{recover}
\wt{\la}_t = \wt{\Cou} \circ \wt{P}_t, \quad t \geq 0.
\end{equation}
Moreover, the semigroup $(\wt{P}_t)_{t\geq 0}$ is completely positive and unital
if and only if $(\lambda_t)_{t\geq 0}$ is a convolution semigroup of states.
We stress the point that, in the noncompact case ($\bialg$ nonunital),
$P_t:=R_{\lambda_t}$ need not leave
$\bialg$ invariant and so \emph{there may be no semigroup on $\bialg$ itself}.
In the next section we shall see that,
under a natural condition on $\bialg$, this obstruction is removed.

The semigroups on $\multib$ which are associated with
convolution semigroups of functionals on $\bialg$ are
characterised in several simple ways.
This is the content of the next result.
The convention on composing strict maps
permits us to dispense with almost all tildes.

\begin{tw} \label{assoc}
Let $(\wt{P}_t)_{t \geq 0}$ be a strict semigroup on $M(\bialg)$ and set
$P_t := \wt{P}_t|_{\bialg}: \bialg \to M(\bialg)$ \tu{(}$t\geq 0$\tu{)}.
Then the  following conditions are equivalent:
\begin{rlist}
\item
$(\wt{P}_t)_{t \geq 0}$ is associated to
a convolution semigroup of functionals on $\bialg$;
\item
$(P_t)_{t \geq 0}$ enjoys the commutativity property:
\[
L_{\mu} P_t  = P_t L_{\mu}, \quad \mu \in \bialg^*, t \geq 0;\]
\item
for each $t\geq 0$, $P_t$ is completely bounded,
$\idb \ot P_t$ is strict as a map from $\bialg \ot \bialg$
to $M(\bialg \ot \bialg)$, and
the following strong invariance condition holds:
\begin{equation}\label{2.3}
 \Com  P_t  = (\idb \ot P_t) \Com;
\end{equation}
\item
$(P_t)_{t \geq 0}$ satisfies the weak invariance condition:
\begin{equation}\label{invar}
P_t  = \big(\idb \ot (\Cou\circ P_t)\big) \Com,
\quad t\geq 0.
\end{equation}
\end{rlist}
\end{tw}

\begin{proof}
The implication (iii)$\Rightarrow$(ii) is immediate.

(ii)$\Rightarrow$(iv):
Assume that (ii) holds and let $t\geq 0$.
Then, since $\Cou \circ L_\mu = \mu$,
\[
\mu \circ \big( \idb \ot (\Cou \circ P_t) \big) \Com =
(\Cou \circ P_t) L_\mu = \Cou \circ (L_\mu P_t) =
\mu \circ P_t, \quad \mu \in \bialg^*,
\]
so (iv) holds.

(iv)$\Rightarrow$(i):
Assume that (iv) holds and define
$\lambda_t = \Cou \circ P_t \in \bialg^*$ ($t\geq 0$). Then, for $s,t\geq 0$,
\begin{align*}
(\lambda_s \ot \lambda_t)\Com
&=
((\Cou \circ P_s) \ot (\Cou \circ P_t)) \Com
\\&=
\big(\Cou\circ P_s\big)\big(\idb \ot (\Cou\circ P_t)\big)\Com
=
\big(\Cou \circ P_s\big)P_t  = \Cou \circ P_{s+t}
=
\lambda_{s+t},
\end{align*}
and
\[
R_{\lambda_t} =
(\idb \ot \lambda_t)\Com =
\big(\idb \ot (\Cou\circ P_t)\big)\Com = P_t.
\]
Thus $(\lambda_t)_{t\geq 0}$ is a convolution semigroup of functionals on $\bialg$ and
$(\wt{P}_t)_{t\geq 0}$ is its associated semigroup on $M(\bialg)$.

(i)$\Rightarrow$(iii):
Assume that $(\wt{P}_t)_{t\geq 0}$ is the semigroup on $M(\bialg)$ associated with a
convolution semigroup of functionals $(\lambda_t)_{t\geq 0}$ on $\bialg$, and let $t\geq 0$.
Then $P_t$ is completely bounded,
$\idb \ot P_t$ equals the composition $(\idb \ot \idb \ot \la_t) (\idb \ot \Com)$
and so is strict, and~\eqref{2.3} holds
by coassociativity:
\[
\Com P_t =
\Com (\idb \ot \la_t) \Com =
(\idb \ot \idb \ot \la_t) \Com\!^{(2)} =
(\idb \ot P_t) \Com.
\]
\end{proof}


\begin{rems}

The above proof yields the following useful
characterisations of the range of the $R$-map:
\begin{align*}
\Ran R &= \{ T \in B(\bialg; M(\bialg)): T \text{ is strict and } L_\mu T = T L_\mu \text{ for all } \mu \in \bialg^* \}
\\ &=
\{ T \in CB(\bialg; M(\bialg)):
(\idb \ot T) \text{ is strict and } \Com T = (\idb \ot T)\Com \}
\\ &=
\{ T \in B(\bialg; M(\bialg)): T = (\idb \ot (\Cou \circ T))\Com \}.
\end{align*}
If $\bialg$ is a non-counital $C^*$-bialgebra then we still have
\begin{align*}
\Ran R &\subset \{ T \in CB(\bialg; M(\bialg)): (\idb \ot T) \text{ is strict and } L_\mu T = T L_\mu \text{ for all } \mu \in \bialg^* \}
\\ &=
\{ T \in CB(\bialg; M(\bialg)):
(\idb \ot T) \text{ is strict and } \Com T = (\idb \ot T)\Com \}.
\end{align*}
The equality and inclusion above follow, in turn, from the facts that the set $\{ (\mu \ot \idb)\, \wt{\ } \, : \mu \in
\bialg^* \}$ separates $M(\bialg \ot \bialg)$, and that $(\mu \ot \idb \ot \la)\Com\!^{(2)}$ is a common expression for $L_\mu  R_\la$ and $R_\la L_\mu$ ($\la, \mu \in \bialg^*$).
\end{rems}

So far we have not considered the question of continuity, in $t$,
for convolution semigroups and their associated semigroups. It is
easily seen that one is norm-continuous if and only the other is, in
which case the generating functional and the generator of the
associated semigroup are related by $\tilded{Z } = \wt{R}_\gamma$
and $\wt{\gamma} = \wt{\Cou} \circ \tilded{Z }$. The corresponding
statement at the level of weak/strong continuity need not be true.
For example the translation semigroup on $C_{\bd}(\br)$ is the
semigroup associated with the convolution semigroup of Dirac
measures $(\delta_t)_{t\geq 0}$ on $\br$. Whereas the latter is
weakly continuous, the former is not. However the translation
semigroup leaves $C_0(\br)$ invariant and restricts to a strongly
continuous semigroup there. This is a \emph{Feller property}
(\cite{Jacob}) which, 
in our framework, corresponds to the semigroup $(\wt{P}_t)_{t\geq 0}$
leaving the $C^*$-algebra $\bialg$ invariant and
restricting to a $C_0$-semigroup on $\bialg$.
In the next section we shall see that this holds for all weakly continuous
semigroups on a wide class of locally compact quantum semigroups
which includes all locally compact (quantum) groups.
We warn the reader that there is a variety of definitions
of `Feller' and `strong Feller'---in both the classical and noncommutative
literature (e.g.\ \cite{Sau}).

\section{Convolution semigroups on $C^*$-bialgebras satisfying \\
the residual vanishing at infinity condition}
\label{rvai}
\emph{Throughout this section we assume} that $\bialg$ is a $C^*$-bialgebra which enjoys
the following `residual vanishing at infinity'  property:
\[
(\bialg \ot 1) \Com(\bialg) \subset \bialg \ot \bialg \ \text{ and } \
 (1 \ot \bialg) \Com(\bialg) \subset \bialg \ot \bialg.
\]
The $C^*$-bialgebras arising from locally compact quantum groups in the sense of Kustermans and Vaes (\cite{KV}) satisfy this
property, as do (trivially) all unital $C^*$-bialgebras. In the classical case of a commutative quantum semigroup (see
Section~\ref{section: cases}), $\bialg = C_0(\semiGp)$ for a locally compact Hausdorff semigroup $\semiGp$, $\bialg \ot \bialg$
and $M(\bialg \ot \bialg)$ are identified with $C_0(H \times H)$ and $C_\bd (H \times H)$ respectively, and the coproduct is
given by $\Com(F)(h_1,h_2)=F(h_1 h_2)$, for $F\in\bialg$ and $h_1, h_2\in H$. To see how the condition applies then, let $F,F'\in
C_0(\semiGp)$. If $|F(h)|, |F'(h)| < \epsilon$ for $h\in\semiGp\setminus K$ then $|F(h_1)F'(h_1h_2)| < \epsilon^2$ for
$(h_1,h_2)\in\semiGp \times\semiGp \setminus \phi^{-1}(K\times K)$ where $\phi$ is the continuous map
$(h_1,h_2)\mapsto(h_1,h_1h_2)$. Thus if $\phi$ is proper then $(F\ot 1)\Com F' \in C_0(\semiGp \times\semiGp ) = \bialg\otimes
\bialg$ and a similar argument applies to $(1\ot F)\Com F'$. Lance credits Iain Raeburn for the suggestive terminology
(in~\cite{Lance}).

The maps $R_\mu$ and $L_{\mu}$ defined
in~\eqref{L and R} now act on $\bialg$, as we show next.

\begin{prop}
Let $\mu\in\bialg^*$. Then
\[
R_\mu(\bialg)\subset \bialg \text{ and }
L_{\mu}(\bialg) \subset \bialg.
\]
\end{prop}
\begin{proof}
This follows from the fact, commonly used in topological quantum group theory,
that $\mu \in \bialg^*$ may be factorised as $\lambda c$,
for some $\lambda\in\bialg^*$ and $c\in\bialg$ (\cite{FactLem});
\begin{align*}
L_{\mu} (a) &= (\mu \ot \idb)\, \wt{\ }\, (\Com(a))
\\ &=
(\lambda \ot \idb)\, \wt{\ }\, ((c \ot 1) \Com(a))
\in (\lambda \ot \idb)(\bialg\ot\bialg)\subset \bialg,
\quad a\in\bialg,
\end{align*}
and similarly for $R_\mu$.
\end{proof}
The significance of this is that
each convolution semigroup of functionals $(\lambda_t)_{t\geq 0}$  on $\bialg$ now has
an \emph{associated semigroup on the $C^*$-algebra} $\bialg$ itself:
$(P_t := R_{\lambda_t})_{t\geq 0}$.
Reconstruction of the convolution semigroup from its associated semigroup
now reads:
\begin{equation} \label{simplifies}
\la_t = \Cou \circ P_t, \quad t \geq 0,
\end{equation}
and Theorem~\ref{assoc} now has the following version.

\begin{prop} \label{rvaiassoc}
Let $(P_t)_{t \geq 0}$ be a strict semigroup on $\bialg$.
Then the following conditions are equivalent:
\begin{rlist}
\item
$(P_t)_{t \geq 0}$
is associated to a convolution semigroup of functionals on $\bialg$\tu{;}
\item
$(P_t)_{t \geq 0}$ enjoys the commutativity property
\[
L_{\mu}  P_t  = P_t L_{\mu}, \quad \mu \in \bialg^*, t \geq 0;\]
\item
for each $t\geq 0$, $P_t$ is completely bounded, $\idb \ot P_t$ is strict
and the following strong invariance condition holds:
\[
 \Com  P_t  = (\idb \ot P_t) \Com;
 \]
 \item
$(P_t)_{t \geq 0}$ satisfies the weak invariance condition:
\[
P_t  = (\idb \ot (\Cou \circ P_t)) \Com, \quad t\geq 0.
\]
\end{rlist}
\end{prop}

Moreover we now have a satisfactory correspondence
between continuity properties of the respective semigroups,
showing that
the semigroup on $M(\bialg)$ associated with any weakly continuous
convolution semigroup of functionals on $\bialg$ is necessarily Feller.

\begin{tw}  \label{posgen}
Let $(P_t)_{t \geq 0}$ be the semigroup on $\bialg$ associated with
a convolution semigroup of functionals $(\la_t)_{t \geq 0}$ on $\bialg$.
Then the following are equivalent:
\begin{rlist}
\item
$(\la_t)_{t \geq 0}$ is weakly continuous;
\item
$(P_t)_{t \geq 0}$ is a $C_0$-semigroup.
\end{rlist}
\end{tw}
\begin{proof}
This follows from the identities
\[
(\la_t - \Cou) = \Cou \circ (P_t - \id_\bialg) \text{ and }
\mu \circ (P_t - \id_\bialg) = (\la_t - \Cou) \circ L_\mu,
\]
for $t\geq 0$ and $\mu \in \bialg^*$.
\end{proof}

Under the weak continuity assumption, Proposition \ref{rvaiassoc}
can be formulated at the level of generators.

\begin{prop}
Let
$(P_t)_{t \geq 0}$
be a $C_0$-semigroup on $\bialg$ with generator $\gen$.
Then the following are equivalent:
\begin{rlist}
\item
$(P_t)_{t \geq 0}$ is the semigroup associated
with a weakly continuous convolution semigroup of functionals;
\item
$\gen$ enjoys the commutativity property:
\[
L_\mu \gen \subset \gen L_\mu, \quad \mu \in \bialg^*,
\]
that is,
$L_\mu(\Dom \gen) \subset \Dom \gen$ and
$L_{\mu} \gen (a) = \gen L_{\mu}(a)$ for all
$a\in\Dom \gen$.
\end{rlist}
\end{prop}
\begin{proof}
This follows from Theorem~\ref{posgen},
the equivalence (i)$\Leftrightarrow$(ii) of Proposition~\ref{rvaiassoc},
and Theorem 1.15 of~\cite{Davies}.
\end{proof}

\begin{deft}
For a weakly continuous convolution semigroup of functionals
$(\la_t)_{t \geq 0}$ on $\bialg$, the functional
$\gamma: \Dom\, \gamma \subset \bialg \to \bc$ defined by
\begin{align*}
&\Dom\, \gamma:=
\big\{ a \in \bialg:
\lim_{t \to 0^+} \frac{\la_t (a) - \Cou(a)}{t} \textrm{ exists} \big\};
\\&
\gamma(d) :=
\lim_{t \to 0^+} \frac{\la_t (d) - \Cou(d)}{t}, \quad d\in \Dom\, \gamma,
\end{align*}
is called the \emph{generating functional} of $(\la_t)_{t \geq 0}$.
\end{deft}

\begin{rems}
If $(\la_t)_{t \geq 0}$ is norm-continuous then
the functional $\gamma$ defined above is equal to that of Proposition~\ref{char}, so
our terminology is consistent.
If each $\lambda_t$ is *-preserving then $\Dom\, \gamma$ is selfadjoint and $\gamma$ is *-preserving.
If each $\lambda_t$ is a state,
then $\gamma$ is conditionally positive, in the sense that
$\gamma(a) \geq 0$ for all $a \in \bialg_+\cap \Dom\, \gamma \cap \Ker\, \Cou$
(again, this is consistent with the earlier definition).
The identity~\eqref{simplifies} implies
the inclusion $\Cou \circ \gen \subset \gamma$,
where $\gen$ is the generator of the associated $C_0$-semigroup on $\bialg$,
so $\Dom\, \gamma$ is norm-dense in $\bialg$.

Note that $\gamma$ can be defined
for a weakly continous convolution semigroup on any $C^*$-bialgebra,
however, without the residual vanishing at infinity property,
there is no guarantee that $\gamma$ be densely defined.
\end{rems}

We now recall a useful construction from semigroup theory.
Let  $(P_t)_{t \geq 0}$  be a $C_0$-semigroup on a Banach space $X$
with generator $\gen$.
Set
\begin{equation}\label{Dom}
\Dom \widehat{\gen} :=
\big\{x \in X:
\forall_{\mu \in X^*}
\lim_{t\to 0^+} \frac{\mu(P_t(x)) - \mu(x)}{t} \textrm{ exists} \big\}.
\end{equation}
Clearly $\Dom \widehat{\gen} \supset \Dom \gen$.
By the Uniform Boundedness Principle,
it follows that, for all $x \in \Dom\widehat{\gen}$,
there is a unique element $l \in X^{**}$ satisfying
\[
\lim_{t\to 0^+} \frac{\mu(P_t(x)) - \mu(x)}{t} = l (\mu),
\quad \mu \in X^*.
\]
Setting $\widehat{\gen}(x):= l$,
we have an operator $\widehat{\gen}: X \to X^{**}$ with dense domain
$\Dom \widehat{\gen}$ (justifying the notation~\eqref{Dom}).
Under the canonical identification of $X$ as a subspace of $X^{**}$,
clearly $\widehat{\gen}\supset \gen$. The precise situation is summarised next.
This is a rewording of Theorem 1.24 of~\cite{Davies},
where the reader is warned of a missing hypothesis (not needed here),
namely that---in the notation used there---$\mathcal{L}$
must also be assumed to be invariant under each $T_t^*$.

\begin{lem} \label{Lhat}
Let $(P_t)_{t \geq 0}$ be a $C_0$-semigroup on a Banach space $X$
with generator $\gen$, and associated operator
$\widehat{\gen}:\Dom \widehat{\gen} \subset X \to X^{**}$.
Then
\begin{equation}
\Dom \gen = \{x \in \Dom \widehat{\gen}: \widehat{\gen} (x) \in X\} \textrm{ and }
\gen = \widehat{\gen}|_{\Dom \gen}.
\label{LL}
\end{equation}
In particular,
the operator $\widehat{\gen}$ uniquely determines the semigroup $(P_t)_{t \geq 0}$.
\end{lem}
%

The significance of the operator $\widehat{\gen}$
for present considerations is indicated by the following result.

\begin{prop} \label{Ldet}
Let $(P_t)_{t \geq 0}$ be the semigroup on $\bialg$ associated with
a weakly continuous convolution semigroup of functionals $(\la_t)_{t \geq 0}$ on $\bialg$.
Then its associated operator
$\widehat{\gen}: \Dom \widehat{\gen}\subset \bialg\to\bialg^{**}$ is given by
\begin{align}
&\Dom \widehat{\gen} =
\{ a \in \bialg: \forall_{\mu\in\bialg^*}\  L_\mu(a) \in \Dom\, \gamma \};
\label{3.3}\\ 
&\widehat{\gen} (a) (\mu) = \gamma(L_{\mu}(a)),
\quad a\in\Dom \widehat{\gen}, \mu\in\bialg^*,
\label{form}
\end{align}
where $\gamma$ is the generating functional of $(\lambda_t)_{t\geq 0}$.
\end{prop}
\begin{proof}
If $\mu \in \bialg^*$, $a \in \bialg$ and $t \geq 0$
then $\mu (P_t(a))= \la_t (L_{\mu}(a))$;
in particular $\mu(a) = \Cou (L_{\mu}(a))$.
The identities~\eqref{3.3} and~\eqref{form} now
follow directly from the definitions of $\gamma$ and $\widehat{\gen}$.
\end{proof}

The next result,
the first part of which fully justifies the term `generating functional',
now follows easily from Proposition~\ref{Ldet},
Lemma~\ref{Lhat} and the identity~\eqref{simplifies}.
\begin{tw}
Let $(\la_t)_{t \geq 0}$ be a weakly continuous convolution semigroup of functionals on $\bialg$
with generating functional $\gamma$.
Then
\begin{alist}
\item
the functional $\gamma$ determines the semigroup $(\la_t)_{t \geq 0}$ uniquely;
\item
$\gamma$ is bounded if and only if
$(\la_t)_{t \geq 0}$ is norm-continuous, in which case
$\gamma$ is everywhere defined.
\end{alist}
\end{tw}
%

\section{Convolution semigroups on $C^*$-bialgebras of discrete type}
\label{disc type}
\emph{For this section we assume that $\bialg$ is of discrete type}, that is,
as a $C^*$-algebra,
$\bialg$ is a $c_0$-direct sum of matrix algebras.
The existence of a counit forces $\bialg$ to have the form
\begin{equation}
\bialg = \bc \Omega \oplus \bialg_0 \text{ where }
\bialg_0 = \bigoplus_{\lambda \in \La} M_{n(\la)},
\label{disc}
\end{equation}
for an index set $\La$ and family $(n(\la))_{\lambda\in\Lambda}$ in $\bn$;
the counit being given by the formula
\[
\Cou\Big(\alpha \Omega \oplus \bigoplus_{\la\in \La} a_{\la}\Big) = \alpha.
\]
We identify $\bialg_0$ with $0 \oplus \bialg_0  \subset \bialg$
and view $\bialg$ as represented on the Hilbert space
$\bc \oplus \bigoplus_{\la \in \La} \bc^{n(\la)}$.
Then $M(\bialg)$ equals the $l^\infty$-direct sum
\[
\mlg := \bc \Omega \oplus \prod_{\la \in \La}M_{n(\la)},
\]
which is a von Neumann algebra,
coinciding with the universal enveloping von Neumann algebra of $\bialg$.

The examples we have in mind are
the algebras of functions vanishing at infinity on a discrete semigroup, and
the duals of compact (quantum) groups.
In the former case $n(\la)=1$ for each $\la \in \La$;
in the latter case $\La$ is the set of
equivalence classes of irreducible representations of the underlying compact group --
more on these in the next section.


\begin{lem}
\label{str}
Let $(\la_t)_{\geq 0}$
be a weakly continuous convolution semigroup of states on $\bialg$.
Then its associated semigroup
$(\wt{P}_t)_{t\geq 0}$ on $\multib$
is strongly continuous.
\end{lem}
\begin{proof}
Since $M(\bialg)$ is the von Neumann algebra $\mlg$,
and by the same token $M(\bialg\ot\bialg)$ is the von Neumann algebra $\mlg\otol\mlg$,
Proposition~\ref{new} gives maps
\[
L_{\nu} := (\nu \otM \idm) \wt{\Com} : M(\bialg)\to M(\bialg),
\quad \nu\in
M(\bialg)^* =  \mlg^*,
\]
satisfying
\[
\wt{\lambda}\circ L_{\nu} =
\nu \circ(\idm\, \otol\, \wt{\lambda})\wt{\Com} =
\nu \circ \wt{R}_\lambda,
\quad \lambda\in\bialg^*=\mlg_*,\ \nu\in M(\bialg)^*=\mlg^*.
\]
In particular,
\begin{equation}
\label{commut}
\nu \circ \big(\wt{P_t} - \id_{M(\bialg)}\big)  =
\big(\wt{\la}_t - \wt{\Cou}\big)\circ L_{\nu},
\quad t \geq 0, \ \nu \in M(\bialg)^*.
\end{equation}
Therefore, by Proposition~\ref{2.2*}, $(\wt{P}_t)_{t \geq 0}$
is weakly continuous, and so also strongly continuous.
\end{proof}

\begin{rem}
Note how the existence of slice maps for not-necessarily-strict functionals
is crucially used in the above proof (via Proposition~\ref{new}).
\end{rem}

\begin{tw} \label{main}
Let $(\la_t)_{t \geq 0}$ be a weakly continuous convolution semigroup of states
on $\bialg$. Then $(\la_t)_{t \geq 0}$ is norm-continuous.
\end{tw}
\begin{proof}
Define a functional
$\tilded{\gamma } : \Dom\, \tilded{\gamma } \subset M(\bialg) \to \bc$ by
\begin{align*}
&\Dom\, \tilded{\gamma } :=
\big\{ m \in \multib:
\lim_{t \to 0^+} \frac{\wt{\la}_t (m) - \wt{\Cou}(m)}{t}
\textrm{ exists} \big\};
\\&
\tilded{\gamma } (d) =
\lim_{t \to 0^+} \frac{\wt{\la}_t (d) - \wt{\Cou}(d)}{t},
\quad d\in \Dom\, \tilded{\gamma }.
\end{align*}
Clearly $\Dom\, \tilded{\gamma }$ is a selfadjoint subspace of $M(\bialg)$ containing $1$
and $\tilded{\gamma }(1)=0$. Moreover, the identity~\eqref{recover} and
Lemma~\ref{str} imply that $\Dom\, \tilded{\gamma }$ is norm-dense in $M(\bialg)$
(cf.\ the remarks in Section~\ref{rvai}).

Let $p$ be the orthogonal projection $1 - \Omega \in M(\bialg)$.
Note that $p$ is the identity of the ideal $\Ker\, \wt{\Cou}
= M(\bialg_0)$ of $M(\bialg)$.
By density and selfadjointness of $\Dom\, \tilded{\gamma }$
there exists  $z \in \Dom\, \tilded{\gamma }$
such that $z=z^*$ and $\|p-z\| \leq \frac{1}{4}$.
Note that $|\wt{\Cou}(z)| = |\wt{\Cou}(z-p)|\leq \frac{1}{4}$.
Put $y = z - \Cou(z) 1$.
Then $y \in \Dom\, \tilded{\gamma }$ and $\wt{\Cou}(y)=0$, so in particular
$\sup_{t>0}\frac{\wt{\la}_t(y)}{t} <\infty$.
Moreover, since $y \in M(\bialg_0)$ and $z\geq p - \frac{1}{4}1$,
\[
y = pyp =
(pzp-\wt{\Cou}(z)p)\geq p\big( p-\frac{1}{4}1\big)p - \frac{1}{4}p = \frac{1}{2} p.
\]
Therefore
\[
C:=\sup_{t>0} \frac{\wt{\la}_t(p)}{t} \leq  2 \sup_{t>0} \frac{\wt{\la}_t(y)}{t} <\infty.
\]
Now, for each $t >0$, define a positive functional $\mu_t$ in $\strfunb$ by the formula
\[
\mu_t (m) = \frac{1}{t} \wt{\la}_t (pmp), \quad m \in M(\bialg).
\]
All these functionals have norm smaller than $C$:
\[
\|\mu_t\| = \mu_t(1) = \frac{1}{t} \wt{\la}_t (p) \leq C.
\]
Now let $\nu$ be an arbitrary functional in $\multib^*$.
Then, since $p$ is the identity element of $M(\bialg_0)$,
it follows from~\eqref{commut} that, for all $m \in M(\bialg)$,
\begin{align*}
\nu\big(\wt{P}_t(m) - m\big)  &=
\wt{\lambda}_t \big(L_{\nu}(m)- \wt{\Cou}(L_{\nu}(m))1\big)
\\ &=
\wt{\lambda}_t \big( p\big( L_{\nu}(m) - \wt{\Cou}(L_{\nu}(m))1\big) p\big)
=
t \,\mu_t \big(  L_{\nu}(m) - \wt{\Cou}(L_{\nu}(m))1\big).
\end{align*}
Thus,
letting $\tilded{Z }$ be the generator of the $C_0$-semigroup $(\wt{P}_t)_{t\geq 0}$,
\begin{align*}
|\nu(\tilded{Z } (m))| = &
\lim_{t \to 0^+} \frac{1}{t}
\left|\nu\big(\wt{P}_t(m) - m\big)\right|
\\ \leq &\
2 C\| L_\nu(m)\|
\leq 2 C\|\nu\| \|m\|,
\quad m \in \Dom \tilded{Z } .
\end{align*}
It follows that $\tilded{Z }$ is bounded, so
$(\wt{P}_t)_{t \geq 0}$ is norm-continuous,
and therefore $(\la_t)_{t \geq 0}$ is too.
\end{proof}

\begin{rem}
Since the multiplier algebra $M(\bialg)$ is a von Neumann algebra,
Theorem~\ref{main} can alternatively be deduced from Lemma~\ref{str}
using the fact that any strongly continuous completely positive semigroup
on a von Neumann algebra is automatically norm-continuous (\cite{George}).
\end{rem}

We note here that the results of
Sections~\ref{on Cstar Bialgebras},~\ref{rvai} and~\ref{disc type}
remain valid for a \emph{$C^*$-hyperbialgebra} $\bialg$,
where the coproduct is only assumed to be completely positive, strict and preunital.
Examples of $C^*$-hyperbialgebras
are provided by compact quantum hypergroups (\cite{ChV}).
The classical theory is described in~\cite{BloomHeyer}.

\section{Commutative and cocommutative cases}
\label{section: cases}
Let $\bialg$ be a commutative $C^*$-bialgebra of discrete type.
Then it follows from Gelfand theory that $\bialg$ is of the form
$C_0(\Gamma)$, where $\Gamma$ is a discrete semigroup,
and the coproduct and counit are given by
\[
\Com(F) (\gamma,\gamma') = f(\gamma \gamma')
\text{ and }
\Cou(F) = F(e), \quad F \in C_0(\Gamma), \gamma, \gamma' \in \Gamma,
\]
under the natural identifications $M(\bialg \ot \bialg) = M(C_0(\Gamma \times \Gamma)) = C_\bd(\Gamma\times \Gamma)$. If $\Gamma$
satisfies both left and right cancellation properties then  the map $\phi$  introduced in the first paragraph of Section
\ref{rvai} is proper, and so is its `right' version, so that the residual vanishing at infinity property holds. The convolution
semigroups of states on $\bialg$ correspond to convolution semigroups of probability measures on $\Gamma$ via the Riesz
Representation Theorem, and Theorem~\ref{main} specialises as follows.

\begin{prop}
Every convolution semigroup of probability measures on a discrete semigroup has bounded infinitesimal generator.
\end{prop}

\noindent
This extends the central conclusion of Theorem 4.1.5 of~\cite{Heyer},
which is formulated for discrete groups.

Recall that \emph{cocommutativity} for a $C^*$-bialgebra $\bialg$ means
$\Sigma\Com = \Com$ where $\Sigma$ denotes the tensor flip
on $\bialg \ot \bialg$, which is strict as a map $\bialg \ot \bialg \to M(\bialg \ot \bialg)$;
it is equivalent to
commutativity of the convolution product on $\bialg^*$.
Now, for a locally compact group $G$,
the reduced and universal group $C^*$-algebras are isomorphic
via the left regular representation if and only if
$G$ is amenable (\cite{Pedersen}, Theorem 7.3.9),
in particular this holds if $G$ is compact.
Therefore, by
the Peter-Weyl theory of unitary representations of compact groups,
the group $C^*$-algebra of a compact group is isomorphic to
a $c_0$-direct sum of matrix algebras.

We are not aware of
any characterisation of cocommutative $C^*$-bialgebras of discrete type.
However every cocommutative discrete quantum group $\bialg$
is the group $C^*$-algebra of a compact group $G$ where, viewing
$C^*(G)$ as the $C^*$-subalgebra of $B(L^2(G, \text{Haar}))$
generated by the set of convolution operators
$\lambda(F): \xi \mapsto F \star \xi$ ($F \in C(G)$),
the coproduct and counit are determined by
\[
\wt{\Com} (\lambda_g) = \lambda_g \ot \lambda_g
\text{ and } \wt{\Cou} (\lambda_g) = 1,
\quad g \in G,
\]
$\lambda_g\in M(\bialg)$ being the translation operator
$(\lambda_g \xi)(g') = \xi(g^{-1}g')$.
This follows from duality theory for locally compact quantum groups
(\cite{BB}, \cite{Kuslect}). Specifically, discrete quantum groups are
naturally isomorphic to their quantum group biduals;
the dual of a cocommutative discrete quantum group
is a commutative compact quantum group; and
the dual of a commutative compact quantum group $C(G)$ is
the group $C^*$-algebra of $G$ with the quantum group structure
defined above.

It is easily verified that
the unital *-algebra spanned by $\{\la_g: g \in G\}$
is strictly dense in the multiplier algebra $M(\bialg)$, and that
the map $g \mapsto \la_g$ is strictly continuous.
The following consequence of Proposition 7.1.9 of~\cite{Pedersen}
is noted here for convenience.
\begin{prop}
Let $\bialg$ be the group $C^*$-algebra of a compact group $G$.
Then
\begin{equation}
\label{corresp}
\wt{\omega}(\la_g) = \phi(g), \quad g \in G,
\end{equation}
describes a one-to-one correspondence between the
positive linear functionals $\omega$ on $\bialg$ and
the continuous positive-definite functions $\phi$ on $G$.
\end{prop}

Given a positive-definite function $\phi$ on a group $G$,
the function $\psi: G\to\mathbb{C}$, $g \mapsto \phi(g) - \phi(e)$
plainly enjoys the following three properties:
it is \emph{conditionally positive definite}, that is,
for all $n\in\mathbb{N}$,
$g_1, \ldots ,g_n \in G$ and
$z_1, \ldots ,z_n \in \bc$
satisfying $\sum_{i=1}^n z_i = 0$,
\[
\sum_{i,j=1}^n \overline{z_i}z_j \psi(g_i^{-1}g_j) \geq 0,
\]
it is \emph{Hermitian}, i.e.
\[
\psi = \psi^* \text{ where }
\psi^*(g) := \overline{\psi(g^{-1})}, \quad g\in G,
\]
and it vanishes at the identity element $e$:
\[
\psi (e) = 0.
\]
The content of the theorem below is that the converse holds
if $G$ is compact and $\psi$ is continuous.
We deduce this from our results.

\begin{tw}[\cite{Guichardet}, Theorem 4.1]
Let $G$ be a compact group and let $\psi$ be a continuous, Hermitian,
conditionally positive-definite function on $G$ vanishing at the identity element $e$.
Then there is a positive-definite function $\phi$ on $G$ such that
\begin{equation}
\label{diff}
\psi(g)= \phi(g) - \phi(e), \quad g \in G,
\end{equation}
in other words there is a constant $c \in \mathbb{R}_+$ such that
$\psi + c$ is positive-definite.
\end{tw}

\begin{proof}
Let $\bialg$ be the discrete quantum group $C^*(G)$, as above,
and let $\mathcal{B}$ be the *-algebraic span of $\{\la_g: g \in
G\}$ in $M(\bialg)$. By the Sch\"onberg correspondence, the continuous function $e^{t\psi}$ is positive-definite for each $t\geq
0$ (\cite{ParthaSch}, Lemma 1.7). Let $(\mu_t)_{t\geq 0}$ be the corresponding family of states on $\bialg$, defined
via~\eqref{corresp}. Then, for all $g \in G$ and $s,t\in \mathbb{R}_+$,
\[
(\wt{\mu}_s \star \wt{\mu}_t) (\la_g) =
\wt{\mu}_s(\la_g) \wt{\mu}_t (\la_g) =
e^{s\psi(g)} e^{t\psi(g)} =
e^{(s+t)\psi(g)} =
\wt{\mu}_{s+t}(\la_g)
\]
and
\[
\wt{\mu}_0(\la_g) = 1 = \wt{\Cou}(\la_g).
\]
By the strict density of $\mathcal{B}$ in $M(\bialg)$,
it follows that
$(\mu_t)_{t\geq 0}$ is a convolution semigroup of states on $\bialg$.
In view of the identity
\[
\mu_t(\lambda(F))  =
\int_G F(g) e^{t\psi(g)}  d\,g,
\quad t \in \mathbb{R}_+, F \in C(G),
\]
it follows from the compactness of $G$ and Lebesgue's Dominated
Convergence Theorem, that $(\mu_t)_{t\geq0}$ is weakly continuous.
Therefore, by Theorem~\ref{main}, it is norm-continuous and so, by
Proposition~\ref{char}, the strict extension of its generating
functional $\gamma$ satisfies
\[
\wt{\gamma}(\la_g) =
\lim_{t \to 0^+}
t^{-1} \big( \wt{\mu}_t(\la_g) - \wt{\Cou}(\la_g) \big) =
\lim_{t \to 0^+}
t^{-1} \big( e^{t\psi(g)} - 1 \big) = \psi(g), \quad g\in G.
\]
By Theorem 7.3 of~\cite{qscc3}, there is a nondegenerate representation
$(\pi,\hil)$ of $\bialg$ and vector $\eta\in\hil$ such that
$\gamma = \omega_\eta \circ (\pi - \iota \circ \Cou)$ where
$\iota$ is the ampliation $\bc \to B(\hil)$ and $\omega_\eta$ denotes
the vector functional $T \mapsto \langle \eta, T \eta \rangle$ on
$B(\hil)$.
Letting $\phi: G \to \bc$ be the positive-definite function
$g \mapsto (\omega_\eta \circ \wt{\pi})(\la_g)$, we have
$\phi(g) - \phi(e) = \wt{\gamma}(\la_g) = \psi(g)$ for all $g\in G$,
and so the proof is complete.
\end{proof}
\begin{rems}
Guichardet's approach is to reduce the proposition to the vanishing of
the first cohomology group of unitary representations of compact groups
and then to appeal to Theorem 15.1 of~\cite{ParthaSch},
which delivers just that.

 In the opposite direction to that presented here, Guichardet's Theorem can be used to give an alternative proof of Theorem 7.2 of \cite{qscc3} in the special case of multiplier $C^*$-bialgebras of the type $C^*(G)$, where $G$ is a compact group.
\end{rems}

\medskip
\noindent

\emph{Acknowledgement.}
We are grateful to Uwe Franz for useful remarks on an earlier draft,
in particular for bringing Guichardet's Theorem to our attention.


\begin{thebibliography}{KaR}




\bibitem[BaS]
{BB} S.\,Baaj and G.\,Skandalis,
Unitaires multiplicatifs commutatifs,
\emph{C.R.\ Math.\ Acad.\ Sci.\,Paris}
\textbf{336}  (2003),  no. 4, 299--304.

\bibitem[Ber]{Bertoin}
J.\,Bertoin,
``L\'evy Processes,''
Cambridge Tracts in Mathematics \textbf{121},
CUP, 1996.

\bibitem[BlH]{BloomHeyer}
W.R.\,Bloom and H.\,Heyer,
``The Harmonic Analysis of Probability Measures on Hypergroups,''
de Gruyter, Berlin, 1995.


\bibitem[ChV]{ChV}
Yu.\,Chapovsky and L.\,Vainerman,
Compact quantum hypergroups,
\emph{J.\ Operator Theory} \textbf{41} (1999) no.\,2, 261--289.

\bibitem[Dav] {Davies}
E.B.\,Davies,
``One-Parameter Semigroups,''
London Mathematical Society Monographs \textbf{15},
Academic Press, London, 1980.

\bibitem[Ell] {George}
G.\,Elliott,
On the convergence of a sequence of completely positive maps to the identity,
\emph{J.\ Austral.\ Math.\ Soc.\ Ser.\ A}
\textbf{68} (2000) no.\,3, 340--348.

\bibitem[Fra]{Uwe}
U.\,Franz,
L\'evy processes on quantum groups and dual groups, \emph{in}
``Quantum Independent Increment Processes, Vol. \!II:
Structure of Quantum L\'evy Processes, Classical Probability and Physics,''
LNM \textbf{1866},
Springer, Heidelberg, 2005.

\bibitem[Gui] {Guichardet}
A.\,Guichardet,
``Symmetric Hilbert Spaces and Related Topics,''
Lecture Notes in Mathematics
\textbf{267}, Springer, Heidelberg, 1970.

\bibitem[Hey] {Heyer}
H.\,Heyer,
``Probability Measures on Locally Compact Groups,''
Springer-Verlag, Berlin, 1977.

\bibitem[Jac] {Jacob}
N.\,Jacob,
``Pseudo-differential Operators and Markov Processes, Vol.\,I:
Fourier Analysis and Semigroups,''
Imperial College Press, London, 2001.

\bibitem[KaR]{KaR}
R.V.\,Kadison and J.R.\,Ringrose,
``Fundamentals of the theory of operator algebras, Vol.\,II:
Advanced Theory,''
Graduate Studies in Mathematics \textbf{16},
AMS,
1997.

\bibitem[Ku$_1$] {Johan}
J.\,Kustermans,
One-parameter representations on $C^*$-algebras,
\emph{Preprint, Odense Universitet}, arXiv: \#funct-an/9707010.

\bibitem[Ku$_2$]{Kuslect}
--- --- ,
Locally compact quantum groups,
\emph{in},
``Quantum Independent Increment Processes,
Vol. \!I: From Classical Probability to Quantum Stochastics,"
\emph{eds. \!U. \!Franz \& M. \!Sch\"urmann},
Lecture Notes in Mathematics \textbf{1865},
Springer, Heidelberg, 2005.

\bibitem[KuV]{KV}
J.\,Kustermans and S.\,Vaes,
Locally compact quantum groups,
\emph{Ann.\ Sci.\ \'Ecole Norm.\ Sup. (4)}
\textbf{33} (2000) no.\,6, 837--934.

\bibitem [Lan] {Lance}
E.C.\,Lance,
``Hilbert $C^*$-modules,''
London Mathematical Society Lecture Notes Series \textbf{210},
Cambridge University Press, Cambridge, 1995.


\bibitem[$\text{LS}_1$]
{qscc2}
J.M.\,Lindsay and A.G.\,Skalski,
Quantum stochastic convolution cocycles II,
\emph{Comm.\ Math.\ Phys.} \textbf{280} (2008) no.\,3, 575--610;
\bibitem[$\text{LS}_2$]
{qscc3}
--- --- ,
Quantum stochastic convolution cocycles III, \emph{Preprint},  arXiv:0905/2410 [math.QA].

\bibitem[LiW] {SteveMartin}
J.M.\,Lindsay and S.J.\,Wills,
Existence of Feller cocycles on a $C^*$-algebra,
\emph{Bull.\ London Math.\ Soc.} \textbf{33} (2001) no.\,5, 613--621.


\bibitem[Neu] {Neufang}
M.\,Neufang,
Amplification of completely bounded operators and Tomiyama's slice maps,
\emph{J.\ Funct.\ Anal.} \textbf{207} (2004), 300--329.

\bibitem[PaS]{ParthaSch}
K.R.\,Parthasarathy and K.\,Schmidt,
``Positive Definite Kernels, Continuous Tensor Products, and
Central Limit Theorems of Probability Theory,''
Lecture Notes in Mathematics \textbf{272}, Springer, Berlin, 1972.

\bibitem[Ped]{Pedersen}
G.K.\, Pedersen,
``$C^*$-algebras and their automorphism groups,''
London Mathematical Society Monographs \textbf{14},
Academic Press, London New York, 1979.



\bibitem[Sau]{Sau}
J.-L.\,Sauvageot,
Strong Feller semigroups on $C\sp \ast$-algebras,
\emph{J.\ Operator Theory} \textbf{42} (1999) no.\,1, 83--102.

\bibitem[Sch]{Schurmann}
M.\,Sch\"urmann,
``White Noise on Bialgebras,''
Lecture Notes in Mathematics \textbf{1544}, Springer, Heidelberg, 1993.


\bibitem[Tay]{FactLem}
D.C.\,Taylor,
The strict topology for double centralizer algebras,
\emph{Trans.\ Amer.\ Math.\ Soc.} \textbf{150} (1970), 633--643.

\bibitem[Tom]{Tomiyama}
J.\,Tomiyama,
Tensor products and approximation problems of $C^*$-algebras,
\emph{Publ.\ RIMS, Kyoto Univ.} \textbf{11} (1975), 163--183.

\end{thebibliography}
\end{document}